\documentclass[12pt]{amsart}
\usepackage{amsmath,amsfonts,amsthm}
\newtheorem{theorem}{Theorem}
\newtheorem{lemma}[theorem]{Lemma}
\newtheorem*{eisenbudthm}{Theorem (Eisenbud \cite{eisenbud})}
\newtheorem*{maintheorem}{Main Theorem}
\newtheorem{conjecture}{Conjecture}
\theoremstyle{remark}
\newtheorem*{remarks}{Remarks}

\advance\oddsidemargin-.2in
\advance\evensidemargin-.2in
\advance\textwidth.4in

\def\pplogo{\vbox{\kern-\headheight\kern -15pt
\halign{##&##\hfil\cr&{
\ppnumber}\cr\rule{0pt}{2.5ex}&\ppdate\cr} }} \makeatletter
\def\ps@firstpage{\ps@empty \def\@oddhead{\hss\pplogo}%
  \let\@evenhead\@oddhead 
}
\def\maketitle{\par
 \begingroup
 \def\thefootnote{\fnsymbol{footnote}}
 \def\@makefnmark{\hbox
 to 0pt{$^{\@thefnmark}$\hss}}
 \if@twocolumn
 \twocolumn[\@maketitle]
 \else \newpage
 \global\@topnum\z@ \@maketitle \fi\thispagestyle{firstpage}\@thanks
 \endgroup
 \setcounter{footnote}{0}
 \let\maketitle\relax
 \let\@maketitle\relax
 \gdef\@thanks{}\gdef\@author{}\gdef\@title{}\let\thanks\relax}
\renewcommand{\thetable}{\@Alph\c@table}
\makeatother

\def\ppnumber{\vbox{\baselineskip14pt
\hbox{DUKE-CGTP-06-02} 
\hbox{NSF-KITP-05-124}
\hbox{math.AG/0611014} 
}}
\def\ppdate{\vbox{\baselineskip14pt \hbox{November, 2006}\hbox{Revised September, 2008}}}
\date{} 

\begin{document}
\title{Threefold flops via matrix factorization}
\author{Carina Curto and David R. Morrison}
\address{Department of Mathematics, Box 90320, Duke University, 
Durham NC 27708}
\address{Current address (Curto): Courant Institute of
Mathematical Sciences, New York University, New York, NY 10012}
\email{curto@courant.nyu.edu}
\address{Current address (Morrison): Department of Mathematics, University
of California, Santa Barbara, CA 93106}
\email{drm@math.ucsb.edu}
\date{}
\maketitle

The structure of birational maps between algebraic varieties becomes
increasingly complicated as the dimension of the varieties increases.
There is no birational geometry to speak of in dimension one: if two
complete algebraic curves are birationally isomorphic then they are biregularly
isomorphic.  In dimension two we encounter the phenomenon of the blowup
of a point, and every birational isomorphism can be factored into a
sequence of blowups and blowdowns.  In dimension three, however, we
first encounter birational maps which are biregular outside of a subvariety
of codimension two (called the {\em center}\/ of the birational map).  When 
the center has a neighborhood with
trivial canonical bundle, the birational map is called a {\em flop}.  The focus
of this paper will be the case of a {\em three-dimensional
simple flop}, in which the
center is an irreducible curve (necessarily a 
smooth rational curve).
One of the motivations for studying this case is a
 theorem of Kawamata \cite{kawamata:crepant}, 
which says that all birational maps
between Calabi--Yau threefolds can be expressed as the composition
of simple flops
 (in fact, of simple flops between nonsingular varieties).

Important examples of simple flops
 were provided by Laufer \cite{[L]}, 
and
three-dimensional
simple flops were studied in general by Reid \cite{pagoda} and by
Pinkham \cite{[P]}.  One fundamental property
is that the center of the flop can be contracted, leaving a
(singular) variety $X$ which is dominated by both of the varieties
involved in the original flop.  $X$ has a hypersurface
singularity, and thus can be locally described as $\{f=0\}$ for
some polynomial $f$ in
which vanishes at the origin.  In fact, as observed by
Koll\'ar and Mori \cite{[Kol]}, the defining polynomial
 can always be put into the
form
\begin{equation}\label{eq1}
 f = x_1^2 + g(x_2,\dots,x_m)
\end{equation}
in appropriate coordinates, and the flop is then induced by the automorphism
$x_1\to-x_1$ of the singular variety $X$.

However, even in dimension three,
most hypersurfaces of the form \eqref{eq1} do not admit 
{\em simple small resolutions}, that is, blowups
in which the singular point is replaced by a nonsingular rational curve,
so the cited result of Koll\'ar and Mori leaves the classification problem 
open.
In 1992, Katz and the second author proved a detailed 
classification theorem for three-dimensional simple flops
\cite{gorenstein-weyl} in terms of a simple invariant of the singularity:
the {\em length}, which is an integer between $1$ and $6$.
Although this classification
 theorem in principle
gives a complete description of three-dimensional
simple flops, the present work began with a
realization that the existing classification theorem did not
specify sufficient detail about a small neighborhood of the center of
a simple flop to 
answer some fundamental questions arising from string theory.

In the language of string theorists, those questions come from 
the study of a physical
model in which a collection of D-branes is made to wrap a rational
curve which is the center of a simple
flop.  The geometry of a neighborhood of
the center dictates
the physics of the model, but not enough information about such 
neighborhoods was available to answer some important physical questions
about this model.
We were led to the
present work by an analogy between this situation and another situation
considered by string theorists: the so-called {\em Landau--Ginzburg models}\/
in string theory.  These models
 also involve a single polynomial (analogous to $f$) and
also have a D-brane theory which is tricky to describe.  In the case of
Landau--Ginzburg models, Kontsevich \cite{kontsevich} proposed (and
a number of physicists investigated, beginning with \cite{Kapustin:2002bi} and
summarized in \cite{Hori:2004zd}) that the D-branes in the theory 
should be described by the category
of matrix factorizations \cite{eisenbud}
associated to the polynomial.  By analogy, we hoped that matrix factorizations
might be useful in studying simple flops, which has proven to be the case.

The idea behind matrix factorizations is quite simple.  Although an irreducible
polynomial $f(x_1,\dots,x_m)$ can never be factored 
(by definition)
in the ring
$K[x_1,\dots,x_m]$, there may well be {\em matrix factorizations}, that is,
equations of the form
\begin{equation}
\Phi\Psi = f I_k,
\end{equation}
where $\Phi$ and $\Psi$ are $k\times k$ matrices with entries in
$K[x_1,\dots,x_m]$.  A familiar example of matrix factorization is 
provided by the construction of Clifford algebras, the elements of
which can be regarded as determining
matrix factorizations for a quadratic polynomial
$f(x_1,\dots,x_m)=x_1^2+\cdots+x_m^2$.  In this paper, we consider
matrix factorizations for rational double points and their deformations.

The classification theorem for three-dimensional simple flops
 \cite{gorenstein-weyl} is based on Reid's
observation \cite{pagoda}
that the general hyperplane section of the contracted variety
$X$ has a rational double
point singularity, and that each small resolution induces a partial resolution
of that rational double point (dominated by the minimal resolution).
Pinkham \cite{[P]} analyzed the deformation theory of such partial resolutions,
which provided the starting point for the classification
theorem of \cite{gorenstein-weyl}.  In the present work, we begin with the
observation that each partial resolution of a rational double point
(dominated by the minimal resolution) has, by the McKay correspondence,
an associated maximal Cohen--Macaulay module which can be described
as a matrix factorization for the equation of the rational double point.
We conjecture that Pinkham's deformation theory for partial resolutions
is actually a deformation theory of the corresponding matrix factorization,
and we prove this conjecture for rational double points of type
$A_{n-1}$ and $D_n$.  

In the course of the proof, we encounter some
``universal flops'' of lengths $1$ and $2$; the analogous universal
flops of lengths $3$, $4$, $5$, and $6$ are still conjectural (and would
involve the $E_6$, $E_7$ and $E_8$ singularities).  The universal flop
of length $1$ is very well known: it is the hypersurface singularity
\begin{equation} \label{length1}
xy-z^2+t^2=0 ,
\end{equation}
which admits two small resolutions by blowing up the ideals 
$(x,z+t)$
or 
$(y,z+t)$,
respectively.  A theorem of Reid \cite{pagoda} 
guarantees that every flop of length $1$ is locally the pullback of the 
universal
one by a map from $X$ to the hypersurface in equation \eqref{length1}; our
theorem about flops of length $2$ is similar.

Recently, Bridgeland \cite{bridgeland} has given a construction of flops
as describing the passage
 from an initial space $Y$ to a moduli space $Y^+$ for certain
 objects in the derived category of coherent sheaves on $Y$, and this
construction and some related ideas 
have been applied in the physical situation described above
\cite{Aspinwall, Aspinwall-Katz}. Our universal flops, and more generally
the description of flops in terms of matrix factorizations, should enable
one to make Bridgeland's construction and its applications much more
explicit.  (In particular, our approach seems closely related to
van den Bergh's alternate proof \cite{VdB} of Bridgeland's theorem.)
We leave this for future investigation.

The outline of the paper is as follows.  In section \ref{sec:1}, we review the
theory of matrix factorizations, and how the McKay correspondence for
rational double points may be interpreted in terms of the existence
of matrix factorizations.  In section \ref{sec:2}, 
we construct a geometric operation,
the {\em Grassmann blowup}, associated to a matrix factorization, and
show that it can be used to recover the partial resolutions of rational
double points directly from the matrix factorization.  In section \ref{sec:3},
we review the basic technique for studying flops by considering the general
hyperplane section of the contracted variety.  In section \ref{sec:4},
we consider the deformation theory for a matrix factorization, treating
the case of $A_{n-1}$ in detail and then conjecturing results for
the other rational double points.  Section \ref{sec:5} 
is devoted to the construction
of the universal flop of length $2$ and some of its properties.  Finally,
in section \ref{sec:6} we prove our conjectures for the case of rational double
points of type $D_n$.  In an appendix, we list the matrix factorizations
for rational double points of types $E_6$, $E_7$, and $E_8$, essentially
drawn from \cite{GSV}.

We would like to thank
P.~Aspinwall, 
T.~Bridgeland, 
D.~Eisenbud, 
Y.~Ito, 
A.~Kapustin, 
S.~Katz, 
H.~Kurke, 
M.~Reid, 
and
B.~Ullrich
for helpful discussions during the course of this work.  The first
author was supported by a National Science Foundation Graduate Research
Fellowship and by National Science Foundation grant DMS-9983320;
the second author was supported
by National Science Foundation grants PHY-9907949 and DMS-0301476, 
the Clay Mathematics Foundation, the John Simon Guggenheim
Memorial Foundation, the Kavli Institute for Theoretical Physics
(Santa Barbara), and
the Mathematical Sciences Research Institute (Berkeley) during the preparation
of this paper.
Any opinions,
findings, and conclusions or recommendations expressed in  this material are
those of the authors and do not necessarily reflect the views of the 
National Science Foundation.

\section{Matrix factorizations and the McKay correspondence} \label{sec:1}

A {\em maximal Cohen--Macaulay module over a Noetherian local ring
$R$} is an $R$-module whose
depth is equal to the Krull dimension
of $R$.  We will primarily be interested in the case in which $R$ is
the localization at the origin of $S/(f)$, where
$S=K[x_1,\dots,x_m]$ is a polynomial ring over a field $K$,
and $f\in S$ vanishes at the origin.
For such a ring $R$, a maximal Cohen--Macaulay $R$-module 
can be lifted to an $S$-module
$M$ which is supported on the hypersurface $X=\{f=0\}$; moreover, $M$ 
will be locally free on the smooth locus of $X$.

Any such $M$ can be presented as the cokernel of a map between two
$S$-modules of the same rank:
\begin{equation}
 S^{\oplus k} \overset{\Psi}{\to} S^{\oplus k} \to M \to 0.
\end{equation}
Since $fM=0$, $fS^{\oplus k} \subset \Psi(S^{\oplus k})$ which 
implies that there is a map
\begin{equation}
S^{\oplus k} \overset{\Phi}{\to} S^{\oplus k}
\end{equation}
such that $\Phi \circ \Psi = f \operatorname{id}$.
The pair $(\Phi,\Psi)$ is called a {\em matrix factorization},
and it turns out that
 $\Psi \circ \Phi = f \operatorname{id}$ as well, so that $(\Psi,\Phi)$
is also a matrix factorization.

\begin{eisenbudthm}
There 
is a one-to-one correspondence between iso\-mor\-phism classes of
maximal Cohen--Macaulay modules over $R$ with no free summands
and matrix factorizations $(\Phi,\Psi)$ with no summand of the
form $(1,f)$ (induced by an equivalence between appropriate categories).
\end{eisenbudthm}

Matrix factorizations have been used to give an explicit realization 
of the McKay correspondence.
Let $G$ be a finite subgroup of $SU(2)$.  The quotient singularity
$\mathbb{C}^2/G$ has a number of interesting properties.  First, if
$\pi: Y \to \mathbb{C}^2/G$ is the minimal resolution of singularities,
then $R^i\pi_*\mathcal{O}_Y=0$ for $i>0$, i.e., the singularity is
{\em rational.}  Moreover, each 
$\mathbb{C}^2/G$ is  a hypersurface singularity of multiplicity $2$---hence
the name {\em rational double point,} since these two properties of
the singularity characterize this class of singular points \cite{Durfee}.

The original observation of McKay \cite{McKay} gave a correspondence
between the nontrivial irreducible representations of $G$ and the components
of $\pi^{-1}(0)$.
The correspondence was made very explicit by 
Gonzalez-Sprinberg and Verdier \cite{GSV} as follows.  The action of
the group $G$
on the coordinate ring $\mathbb{C}[s,t]$ of $\mathbb{C}^2$ makes
$\mathbb{C}[s,t]$ into a $\mathbb{C}[s,t]^G$-module.  As such, it can
be decomposed as
\begin{equation}
 \mathbb{C}[s,t] = \bigoplus_\rho \mathbb{C}[s,t]_\rho ,
\end{equation}
where $\rho:G\to\operatorname{GL}(V_\rho)$ 
runs over the irreducible (complex) representations of $G$.  Each
summand $\mathbb{C}[s,t]_\rho$ takes the form 
$M_\rho \otimes_{\mathbb{C}}V_\rho$, where $M_\rho$ is a maximal 
Cohen--Macaulay module for
the ring $R=\mathbb{C}[s,t]^G$.
Gonzales-Sprinberg and Verdier computed
an explicit presentation of each such module, giving a 
matrix factorization\footnote{The fact that the computation gives a matrix
factorization was not explicitly pointed out in the original paper, but became
an important ingredient in subsequent work of Kn\"orrer \cite{CMI}
and Buchweitz--Greuel--Schreyer \cite{CMII}.}
associated to $M_\rho$.
They went on to show that, while $M_\rho$ is not itself locally free,
the pullback $\pi^*(M_\rho)/(\operatorname{tors})$ is locally 
free, and for each
component $E_i$ of $\pi^{-1}(0)$, 
\begin{equation}
E_i\cdot c_1(\pi^*(M_\rho)/(\operatorname{tors})) = \delta_{ij_\rho},
\end{equation}
where $E_{j_\rho}$ is the component of $\pi^{-1}(0)$ associated to $\rho$
by the McKay correspondence.

To take a simple example, if
 $G=\mathbb{Z}/n\mathbb{Z}\subset SU(2)$, the ring $\mathbb{C}[s,t]^G$
is generated by $x=s^n$, $y=t^n$ and $z=st$ subject to the relation
$xy-z^n=0$. (This shows that the quotient is a hypersurface singularity.)
We label the representations of $G$ so that $\rho_k$ acts 
on $V_k \cong \mathbb{C}$ as $u\mapsto
e^{2\pi ik/n}u$.  It is not hard to see that
$M_{\rho_k}:=\mathbb{C}[s,t]_{\rho_k}$ is generated 
over $\mathbb{C}[s,t]^G$ by $s^k$
and $t^{n-k}$, with relations
\begin{equation}
\begin{split}
y s^k &= z^k t^{n-k}\\
z^{n-k} s^k &= x t^{n-k} .
\end{split}
\end{equation}
This leads to a presentation of $M_{\rho_k}$ as the cokernel of the matrix
\begin{equation}
\begin{bmatrix}
y & -z^k \\ -z^{n-k} & x
\end{bmatrix} .
\end{equation}
We can easily complete this to a matrix factorization:
\begin{equation} \label{eq:an-mat}
\begin{bmatrix}
x & z^{k} \\ z^{n-k} & y
\end{bmatrix}
\begin{bmatrix}
y & -z^k \\ -z^{n-k} & x
\end{bmatrix}
=(xy-z^n)I_2 .
\end{equation}

In general, the Gonzalez-Sprinberg--Verdier matrix factorizations have
the following property: if the locally free module 
$\pi^*(M_\rho)/(\operatorname{tors})$
has rank $\ell_\rho$, then the matrix factorization is given by a pair of
$2\ell_\rho\times2\ell_\rho$ matrices. Moreover, the integers
$\ell_\rho$ are determined by the McKay correspondence and the relation
\begin{equation}
\pi^{-1}(\mathfrak{m}) = \sum \ell_\rho E_{j_\rho},
\end{equation}
where $\mathfrak{m}$ is the maximal ideal of $0\in \mathbb{C}^2/G$.

The Gonzalez-Sprinberg--Verdier matrix factorizations of type $A_{n-1}$
were given above, and those of type $D_n$ can be obtained from the
matrix factorizations of section~\ref{sec:6} by specialization of
parameters (see equations \eqref{eq:matfacDneven} and \eqref{eq:matfacDnodd}
below). 
 For completeness, we have listed the Gonzalez-Sprinberg--Verdier
matrix factorizations for $E_6$, $E_7$, and $E_8$ in an appendix to this
paper.

\section{Grassmann blowups} \label{sec:2}

We now describe a blowup associated to a matrix factorization (or a more
general resolution of sheaves).  Given a
matrix factorization $(\Phi,\Psi)$, the cokernel of $\Psi$ is supported
on the hypersurface $f=0$; since we are assuming our hypersurface $f=0$ is
reduced and irreducible, the  $k\times k$ matrix
$\Psi$ will have some generic
rank $k-r$ along the hypersurface, and the cokernel of $\Psi$ will have
rank $r$.

More generally, we can consider maps $\Psi$ whose cokernel has proper support
$X\subset \mathbb{C}^n$,
such that the cokernel has rank $r$ at each generic point of the support.

We do a {\em Grassmann blowup\footnote{Note that when the module in
question is the pushforward of the tangent bundle from the smooth locus
of $X$, this construction coincides with the more familiar {\em Nash
blowup}.} of $\Psi$}\/
 on which 
there is a locally free sheaf that agrees with $\operatorname{coker} \Psi$
at generic points,
as follows.  In the product
$\mathbb{C}^n \times \operatorname{Gr}(k-r,k)$ we take the closure of
the set\footnote{We could have equally well formulated this definition
in terms of $\operatorname{Gr}(r,k)$ and $\operatorname{coker} \Psi_x$,
but the equivalent formulation we give is more convenient for computation.}
\begin{equation}
 \{ (x,v) \ | \ x \in X_{\text{smooth}},  
v = \operatorname{ker} \Psi_x \}.
\end{equation}
There are natural coordinate charts for this blowup, given by Pl\"ucker
coordinates for the Grassmannian, and a locally free sheaf defined by
pulling back the universal quotient bundle from $\operatorname{Gr}(k-r,k)$.

The following lemma is immediate.

\begin{lemma}
If $\nu:\widetilde{X}\to X$ is any birational map such that 
$\nu^*(\operatorname{coker} \Psi)/(\operatorname{tors})$ is locally free, then
$\nu$ factors through the normalization of the Grassmann blowup of $\Psi$.
\end{lemma}

From this, and the computations in \cite{GSV}, we deduce:

\begin{theorem}
Let $G\subset SU(2)$ be a finite group, let
$\rho$
be a nontrivial
irreducible representation of $G$, let $M_\rho$ be the associated
maximal Cohen--Macaulay $\mathbb{C}[s,t]^G$-module,
and let $\pi_\rho:X_\rho \to \mathbb{C}^2/G$ be the 
(normalization\footnote{In fact, the Grassmann blowup is normal in this case.}
 of the) Grassmann blowup associated to a matrix factorization of $M_\rho$.
Then $X_\rho$ has only rational double point singularities, $\pi_\rho^{-1}(0)$
is an irreducible curve $C_\rho$, and the pullback of $C_\rho$ to the
minimal resolution coincides with $E_{j_\rho}$.
\end{theorem}

In other words, the Grassmann blowup of the 
Gonzalez-Sprinberg--Verdier
matrix factorization directly
produces
the associated curve for the McKay correspondence.

\begin{proof}
Let $\pi: Y \to \mathbb{C}^2/G$ be the minimal resolution
of singularities.
As remarked above, 
Gonzalez-Sprinberg and Verdier 
showed that 
$\pi^*(M_\rho)/(\operatorname{tors})$ is locally free.  Moreover, by rather
intensive computation they showed that
there is a curve $\Gamma_\rho$ meeting $E_{j_\rho}$ in a single point
(and disjoint from other components of $\pi^{-1}(0)$) 
with $c_1(\pi^*(M_\rho))$ supported
on $\Gamma_\rho$.  Thus, for a sufficiently small open neighborhood $U$
of $\overline{\pi^{-1}(0)-E_{j_\rho}}$,
the locally free sheaf
$\pi^*(M_\rho)/(\operatorname{tors})$ 
restricts to a
trivial bundle of rank $\ell_\rho$
on $U$.

It follows that if  $\nu_\rho: \widetilde{X}_\rho
\to \mathbb{C}^2/G$ is the partial resolution
which blows up precisely the curve $E_{i_\rho}$, then
$\pi^*(M_\rho)/(\operatorname{tors})$ pushes forward to a locally free sheaf
on $\widetilde{X}_\rho$, which must coincide with
$\nu_\rho^*(M_\rho)/(\operatorname{tors})$.  Thus, by the lemma, the partial
resolution $\nu_\rho: \widetilde{X}_\rho
\to \mathbb{C}^2/G$ factors through a map $\widetilde{X}_\rho \to X_\rho$
to the normalization of the Grassmann blowup.  Note that 
$X_\rho \to \mathbb{C}^2/G$ is nontrivial since $M_\rho$ is not itself
locally free; since in addition $\nu_\rho: \widetilde{X}_\rho
\to \mathbb{C}^2/G$ has irreducible exceptional set, it follows that the
map $\widetilde{X}_\rho \to X_\rho$ must be an isomorphism,
with $E_{j_\rho}$ mapping to the exceptional set $C_\rho$.
\end{proof}

The proof we have given is unfortunately rather computationally intensive,
insofar as it relies on the detailed explicit computations of \cite{GSV}.
We strongly suspect that a more direct proof should be possible, using the
ideas of Ito--Nakamura \cite{Ito-Nakamura}, but we have not carried
this out.

In the remainder of this paper, we will study deformations of the
Gonzalez-Sprinberg--Verdier matrix factorizations, and their Grassmann
blowups.

\section{Simultaneous resolution and flops} \label{sec:3}

A {\em simple flop} 
is a birational map 
$Y \dashrightarrow Y^+$ 
between Gorenstein threefolds
 which induces an isomorphism
$(Y - C) \cong (Y^+ - C^+),$
 where $C$ and $C^+$ are smooth
rational curves on 
$Y$ and $Y^+$, 
respectively, and 
\begin{equation}
K_Y \cdot C  = K_{Y^+} \cdot C^+ = 0 .
\end{equation}
The curves $C$ and $C^+$ can be contracted
to points (in $Y$ and $Y^+$, respectively), yielding the same normal
variety $X$.

Recall that a
 {\em two-dimensional ordinary double point}\/ is a singular hypersurface
with defining polynomial
\begin{equation}
xy-z^2 .
\end{equation}
Such a singularity has a versal deformation with defining polynomial
\begin{equation} \label{eq:versalA1}
xy-z^2+s=0,
\end{equation}
where $s$ is the coordinate in the versal deformation space.
In 1958, Atiyah \cite{Atiyah} 
noticed that if he made the basechange $s=t^2$ in this
versal deformation,
then the resulting family of surfaces had a {\em simultaneous
resolution of singularities}. 
That is, factoring the equation 
\begin{equation} xy-z^2+t^2=0 \end{equation}
as
\begin{equation}
xy=(z+t)(z-t)
\end{equation}
and blowing up the (non-Cartier) divisor described by 
$x=z+t=0$,
one obtains a family of non-singular surfaces which, for each value
of $t$, resolves the corresponding singular surface.
In fact, there are two ways of doing this, for one might have chosen
to blow up the divisor 
$y=z+t=0$
instead.  This produces two threefolds
$Y$ and $Y^+$ related by a simple flop.  

Reid \cite{pagoda} studied
a generalization of Atiyah's flop with equation of the form
\begin{equation} \label{eq:reidflop}
xy=z^2-t^{2n}=(z+t^n)(z-t^n) ,
\end{equation}
where again the blowups are given by $x=z+t^n=0$ or $y=z+t^n=0$.
Other, more complicated, 
examples of simple flops were given by Laufer \cite{[L]}, and generalized
by Pinkham and the second author \cite{[P]} and by Reid \cite{pagoda}.

A generalization of Atiyah's observation about simultaneous resolutions
was made by Brieskorn \cite{[Bri0]} and Tyurina \cite{[T]},
who showed that for any rational double point, the universal family
over the versal deformation space admits
a simultaneous resolution after basechange.
More precisely, each rational double point has an associated Dynkin 
diagram $\Gamma$ whose Weyl group $\mathcal{W}=\mathcal{W}(\Gamma)$ 
acts on the complexification
$\mathfrak{h}_\mathbb{C}$ of the Cartan subalgebra $\mathfrak{h}$ of the
associated Lie algebra $\mathfrak{g}=\mathfrak{g}(\Gamma)$.  A model
for the versal deformation space is given by 
\begin{equation}
\operatorname{Def} = \mathfrak{h}_\mathbb{C}/\mathcal{W} ,
\end{equation}
and there is a universal family $\mathcal{X}\to \operatorname{Def}$ of 
deformations of the rational double point over that space.

The deformations of the resolution are given by a representable functor,
which can be modeled by
\begin{equation}
\operatorname{Res} = \mathfrak{h}_\mathbb{C},
\end{equation}
and the basechange which relates $\operatorname{Res}$ to $\operatorname{Def}$
coincides with the quotient map 
$\mathfrak{h}_\mathbb{C} \to \mathfrak{h}_\mathbb{C}/\mathcal{W}$.
In fact, there is a universal simultaneous resolution $\widehat{\mathcal{X}}$
of the family $\mathcal{X}\times_{\operatorname{Def}}\operatorname{Res}$.
The construction of this resolution requires some trickiness with the
algebra which will in fact be generalized and somewhat explained later in
this paper.

Pinkham \cite{[P]} adapted 
the work of Brieskorn and Tyurina to cases
in which one does not wish to fully resolve the rational double point, but
only to partially resolve it.  
If we let $\Gamma_0\subset\Gamma$ be the
subdiagram for the part of the singularity that is not being resolved,
then we can define a functor of deformations of the partial resolution,
which has a model
\begin{equation}
\operatorname{PRes}(\Gamma_0) = \mathfrak{h}_\mathbb{C}/\mathcal{W}(\Gamma_0) ,
\end{equation}
and there is a simultaneous partial resolution 
$\widehat{\mathcal{X}}(\Gamma_0)$ of the family
$\mathcal{X}\times_{\operatorname{Def}} \operatorname{PRes}(\Gamma_0)$.

By a lemma of Reid \cite{pagoda}, 
given a simple flop from $Y$ to $Y^+$ and the 
associated small contraction $Y\to X$, the general hyperplane section
of $X$ through the singular point $P$ has a rational double point at $P$,
and the proper transform of that surface on $Y$ gives a partial resolution
of the rational double point (dominated by the minimal resolution).
Pinkham \cite{[P]} used this to give a construction
for all Gorenstein threefold singularities with small resolutions
(with irreducible exceptional set): they
can be described as pullbacks of the universal family
via a map from the disk to 
$\operatorname{PRes}(\Gamma_0)$ 
(for some $\Gamma_0\subset\Gamma$ which is the complement of a
single vertex).

Note that in the examples of Atiyah and Reid, 
one starts with the versal deformation
of $A_1$ given by equation \eqref{eq:versalA1} whose deformation space
$\operatorname{Def}$ has coordinate $s$;  pulling this back via the map
$s=t^{2n}$ yields equation \eqref{eq:reidflop}.
The fact that the power of $t$ is even means that the map $s=t^{2n}$
factors through the degree two cover 
$\operatorname{Res}\to\operatorname{Def}$, as expected from the general
theory.

Koll\'ar \cite{[CKM]}
introduced an invariant of simple flops called the {\em length}:
it is defined to be the generic rank of the sheaf on $C$ defined as the
cokernel of $f^*(\mathfrak{m}_P)\to \mathcal{O}_Y$ where 
$\mathfrak{m}_P$ is the maximal ideal of the singular point $P$.
It is easy to see that the length can be computed from the hyperplane section,
and it coincides with the coefficient of the corresponding vertex in the
Dynkin diagram in the linear combination of vertices which yields the
longest positive root in the root system.
In the Atiyah and Reid cases, the length is $1$.

As mentioned in the introduction, 
Katz and the second author \cite{gorenstein-weyl} classified simple
flops by showing that the 
generic hyperplane section for
a simple flop of length $\ell$ is the smallest rational double point which
uses $\ell$ as a coefficient in the maximal root.  The proof was
computationally intensive, and Kawamata \cite{kawamata:flops}
later gave a short and direct
proof of this result.

\section{Deformations} \label{sec:4}

The versal deformation of the $A_{n-1}$ singularity has defining polynomial
\begin{equation}
 xy-f_n(z) = 0,
\end{equation}
where $f$ is a general monic polynomial 
 of degree $n$ whose coefficient of $z^{n-1}$ 
vanishes.  The coefficients of $f$ are the coordinates on 
$\operatorname{Def}$; the roots of $f$ give coordinates on 
$\operatorname{Res}$ (subject to the constraint that the sum of the
roots is zero).  Note that the action 
of the Weyl group $\mathfrak W_{A_{n-1}}$ coincides with the standard
action of 
symmetric group $\mathfrak S_n$ on the $n$ roots of $f$; the invariants
of this action---the elementary symmetric functions---are 
the coefficients of $f$.

The partial resolution corresponding to the $k^{\text{th}}$ vertex in
the Dynkin diagram corresponds in the invariant theory to the subgroup
$\mathcal{W}(\Gamma_0)\subset \mathcal{W}(\Gamma)$, which in this case is
\begin{equation}
\mathfrak S_k \times \mathfrak S_{n-k} \subset \mathfrak S_n.
\end{equation}
The relationship between the invariants of these two groups is neatly
summarized by writing
\begin{equation}\label{eq:neatsummary}
f_n(z)=g_k(z)h_{n-k}(z),
\end{equation}
where $g$ and $h$ are monic polynomials of degrees $k$ and $n-k$
whose coefficients of $z^{k-1}$
and $z^{n-k-1}$, respectively, sum to zero.  
More precisely, the coefficients
of $f$ give generators for the invariant theory of $\mathcal{W}(\Gamma)$ 
while the
coefficients of $g$ and $h$ give generators for the invariant theory
of $\mathcal{W}(\Gamma_0)$, with the relationship between them specified by
equation \eqref{eq:neatsummary}.
Thus, the coefficients of $g$ and $h$ give coordinates on the partial
resolution space 
$\operatorname{PRes}=\operatorname{Res}/\mathcal{W}(\Gamma_0)$.

It is then easy to see that the matrix factorization \eqref{eq:an-mat}
associated to the $k^{\text{th}}$ vertex of $A_{n-1}$ deforms to
a matrix factorization defined over $\operatorname{PRes}$: just use
\begin{equation}
\begin{bmatrix}
x & g_k(z) \\ h_{n-k}(z) & y
\end{bmatrix}
\begin{bmatrix}
y & -g_k(z) \\ -h_{n-k}(z) & x
\end{bmatrix}
		       =\big(xy-f_n(z)\big)I_2 .
\end{equation}
This matrix factorization encodes the special form of the equation
which was needed in order to find non-Cartier divisors to blow up.
In fact, we can reinterpret the traditional description of these
blowups as being the Grassmann
blowups associated to $\operatorname{coker}(\Psi)$ and
$\operatorname{coker}(\Phi)$, which are normal.
(This reinterpretation of the traditional description of flops of length $1$
is what we will generalize in this paper.)

To see this, we explicitly carry out the Grassmann blowup associated to
the matrix
\begin{equation}
\Psi=
\begin{bmatrix}
y & -g_k(z) \\ -h_{n-k}(z) & x
\end{bmatrix}.
\end{equation}
There are two coordinate charts, in which we use as a basis for 
the kernel of $\Psi$ the vectors $[\alpha \ 1]^T$ and $[1 \ \beta]^T$,
respectively.  In the first coordinate chart, the equation (in matrix
form) defining the variety is
\begin{equation}
\Psi \begin{bmatrix}\alpha\\1\end{bmatrix}
=
\begin{bmatrix}0\\0\end{bmatrix},
\end{equation}
or in other words,
\begin{equation}
\begin{split}
 y\alpha  -g_k(z) &= 0 \\
  -  h_{n-k}(z)\alpha +x &= 0.
\end{split}
\end{equation}
Using the second equation, $x$ can be eliminated, leaving the single equation
\begin{equation}
 y\alpha - g_k(z) = 0 ,
\end{equation}
which defines the versal deformation of an $A_{k-1}$ singularity.

In the second chart, we have the matrix equation
\begin{equation}
\Psi \begin{bmatrix}1\\ \beta\end{bmatrix}
=
\begin{bmatrix}0\\0\end{bmatrix},
\end{equation}
or in other words,
\begin{equation}
\begin{split}
 y - g_k(z)\beta &= 0 \\
  - h_{n-k}(z) + x\beta &= 0.
\end{split}
\end{equation}
This time, the first equation allows us to eliminate $y$, leaving the
single equation
\begin{equation}
 x \beta - h_{n-k}(z) = 0 ,
\end{equation}
which defines the versal deformation of an $A_{n-k-1}$ singularity.

Thus, the Grassmann blowup has produced an exceptional $\mathbb{P}^1$
at two points of
 which the pulled back family has two residual singularities along
the central fiber ($A_{k-1}$ and $A_{n-k}$), and due to the versal deformations
of those singularities appearing on the universal family, every such
deformation can be obtained in this way.

Inspired by this case, and the $D_n$ case to be discussed in 
section~\ref{sec:6}, we formulate the following conjectures:

\begin{conjecture} \label{conj:1}
  For every flop of length $\ell$, there are
two maximal Cohen--Macaulay modules $M$ and $M^+$ 
on 
$X$
of rank $\ell$, 
such that 
$Y$ (resp.\ $Y^+$) is the Grassmann blowup of $M$ 
(resp.\ $M^+$).
\end{conjecture}

\begin{conjecture} \label{conj:2}  For a flop of length $\ell$, the
matrix factorizations corresponding to $M$ and $M^+$ are of size
$2\ell\times2\ell$, and are obtained from each other by switching
the factors $(\Phi,\Psi)\to(\Psi,\Phi)$.
Moreover, if coordinates are chosen so that 
the equation of the hypersurface $X$ takes the
form
$ x_1^2+f(x_2,x_3,x_4)=0$,
then the matrices $\Phi$ and $\Psi$ can be chosen to take the form
\begin{equation}
\Phi = x_1I_{2\ell} - \Xi, \quad \Psi = x_1I_{2\ell} + \Xi ,
\end{equation}
where $\Xi$ is a $2\ell\times2\ell$ matrix whose entries are functions
of $x_2$, $x_3$, $x_4$.  (This is expected thanks to the result of
Koll\'ar and Mori \cite{[Kol]}, which implies that the flop will be
induced by the automorphism $x_1\to-x_1$.)
\end{conjecture}

\begin{conjecture} \label{conj:3} For a partial resolution of a rational
double point corresponding to a single vertex in the Dynkin diagram
with coefficient $\ell$ in the maximal root, 
the versal deformation $\mathcal{X}$ over $\operatorname{PRes}$ has
two matrix factorizations of size $2\ell\times2\ell$, such that the two
simultaneous partial resolutions can be obtained as the Grassmann blowups
of the corresponding Cohen--Macaulay modules.
(These should be regarded as deformations of the Gonzalez-Sprinberg--Verdier
factorizations.)
Moreover, the matrices take the special form $\Phi=xI_{2\ell}-\Xi$,
$\Phi=xI_{2\ell}+\Xi$ in appropriate coordinates.
\end{conjecture}

\begin{maintheorem}  Conjectures \ref{conj:1} and \ref{conj:2} 
hold for lengths $1$ and $2$;
conjecture \ref{conj:3} holds for Dynkin diagrams of type $A_{n-1}$ and $D_n$.
\end{maintheorem}

\begin{remarks} \quad
\begin{enumerate}
\item By the classification theorem of \cite{gorenstein-weyl}, 
in order to prove 
conjectures \ref{conj:1} and \ref{conj:2} for lengths $1$ and $2$,
it suffices to prove
conjecture \ref{conj:3} 
for Dynkin diagrams of type $A_1$ and $D_4$.
\item
It would follow from conjecture \ref{conj:3} 
combined with the classification theorem of
\cite{gorenstein-weyl} that there exists a universal flop of length
$\ell$ for each $1\le\ell\le6$.  Each such universal flop would be obtained
from the universal matrix factorization for Dynkin diagram of length $\ell$
and types $A_1$, $D_4$, $E_6$, $E_7$, $E_8$ (length $5$) and 
$E_8$ (length $6$), respectively.
\item 
Note that we have already given
most of
the proof for the $A_{n-1}$ case of conjecture \ref{conj:3} 
in our discussion above.
We need one additional detail to complete the proof: 
making the change of coordinates $x=u-v$, $y=u+v$
allows us to write the $A_{n-1}$ matrix factorization in the form
$\Phi=uI_2-\Xi$, $\Psi=uI_2+\Xi$, where
\begin{equation}
\Xi=
\begin{bmatrix}
v & h_{n-k}(z) \\ g_k(z) & -v
\end{bmatrix},
\end{equation}
verifying conjecture \ref{conj:3} in this case.
\end{enumerate}
\end{remarks}

\section{The universal flop of length $2$} \label{sec:5}

In this section, we investigate a certain
flop of length $2$ which will turn out to be universal
in an appropriate sense.  Our description of this flop follows an idea of Reid 
\cite[Lemma (5.16)]{pagoda}
although this was not the way we originally found the flop.

We start with a quadratic equation in four variables $x$, $y$, $z$, $t$
over the field $\mathbb{C}(u,v,w)$, chosen so that its discriminant
is a perfect square.  The one we use can be written in  matrix form as
\begin{equation}
\label{eq:quad}
\begin{aligned}
W(x,y,z,t,u,v,w) &:= 
x^2+uy^2+2vyz+wz^2+(uw-v^2)t^2\\
&=\begin{bmatrix}
x&y&z&t\\
\end{bmatrix}
\begin{bmatrix}
1 & 0 & 0 & 0 \\
0 & u & v & 0 \\
0 & v & w & 0 \\
0 & 0 & 0 & uw{-}v^2 \\
\end{bmatrix}
\begin{bmatrix}
x\\y\\z\\t\\
\end{bmatrix},
\end{aligned}
\end{equation}
and its discriminant (the determinant of the matrix of coefficients) is
\begin{equation}
(uw-v^2)^2.
\end{equation}
  (Reid's construction was similar, but had $v=0$;
 as Reid observed, this construction includes both the original Laufer
examples \cite{[L]}
and their generalizations by Pinkham and the second author \cite{[P]}.)  
The general quadratic hypersurface
in $\mathbb{C}^4$ has two rulings by $\mathbb{C}^2$; since
the discriminant is a perfect square, the individual rulings are already
defined over $\mathbb{C}(u,v,w)$.  In fact, the corresponding 
rank two sheaf is a maximal Cohen--Macaulay module over the hypersurface
defined by equation \eqref{eq:quad} (which we now regard as a hypersurface
in $\mathbb{C}^7$).

By Eisenbud's theorem, 
this maximal Cohen--Macaulay
module can be expressed in terms of a matrix factorization $(\Phi,\Psi)$.
We have computed one such, which takes the form
$\Phi = xI_4-\Xi$, $\Psi = xI_4+\Xi$, where
\begin{equation} \label{eq:mainfactor}
\Xi =
\begin{bmatrix}
-vt & y & z & t \\
-uy - 2v z & v t & -ut & z \\
-wz & wt & -vt & -y \\
-uwt & -wz & uy+2v z& v t \\
\end{bmatrix}.
\end{equation}

An explicit computation shows that
\begin{equation}
 -\Xi^2= \left(uy^2 + 2vyz + wz^2 + (uw-v^2)t^2\right)I_4,
\end{equation}
and hence
\begin{equation}
\Phi\Psi =  \left(x^2 + uy^2 + 2vyz + wz^2 + (uw-v^2)t^2\right)I_4.
\end{equation}

A few comments about the geometry are in order.  The quadratic 
form in $x$, $y$,
$z$, $t$ has rank $4$ for generic values of the parameters, but when
$uw-v^2=0$ the rank drops to $2$.  We get rank $1$ at $u=v=w=0$.  This
implies that over each point of $uw-v^2=0$ the fiber of the
Grassmann blowup will be contained in
a $\mathbb{P}^2$ (corresponding to the choice of $\mathbb{C}^2$ within
$\mathbb{C}^3$) and this remains true at $u=v=w=0$ as well.  Moreover,
over $x=y=z=t=0$ we get a fiber contained in $\operatorname{Gr}(2,4)$.
As we will see in detail below, the latter fibers are actually conics
embedded into $\operatorname{Gr}(2,4)$ and all fibers of the Grassmann
blowup have dimension $0$ or $1$.

\begin{theorem}\label{thm:universal}  
For every threefold flop of length $2$, there is a map
from the singular space $X$ to the universal flop of length $2$ such
that the two blowups $Y$ and $Y^+$ are the pullbacks of the Grassmann
blowups of the universal matrices $\Phi$ and $\Psi$.
\end{theorem}

The first step in proving theorem \ref{thm:universal}
is to explicitly carry out the
Grassmann blowup of $\Psi$.
There are six coordinate charts for this blowup, corresponding to
Pl\"ucker coordinates on the Grassmannian, but we will only display
the results in detail for two of these charts.
For the first, we introduce four new variables $\alpha_{i,j}$ and
eight equations $\lambda_{i,j}$ for the blowup by means of 
\begin{equation}
\left[ \lambda_{i,j} \right] = \Psi
\begin{bmatrix}
1 & 0 \\ 0 & 1 \\ \alpha_{1,1} & \alpha_{1,2} \\ \alpha_{2,1} & \alpha_{2,2} \\
\end{bmatrix} =0.
\end{equation}
Note that by multiplying by $\Phi$ from the left, we see that
the original equation is contained in the ideal generated by the
$\lambda_{i,j}$.

It is not difficult to find some elements of the ideal that
are divisible by $z^2+ut^2=\partial W/\partial w$:
\begin{equation}
\begin{aligned}
z\lambda_{1,1}-ut\lambda_{1,2}-t\lambda_{2,1}-z\lambda_{2,2}
&=(z^2+ut^2)(\alpha_{1,1}-\alpha_{2,2}) ,\\
-ut\lambda_{1,1}-uz\lambda_{1,2}-z\lambda_{2,1}+ut\lambda_{2,2}
&=(z^2+ut^2)(2v-\alpha_{1,2}u-\alpha_{2,1}) ,\\
(\alpha_{1,2}uz+\alpha_{2,2}ut)\lambda_{1,2}
+ (-\alpha_{1,2}ut+&\alpha_{2,2}z)\lambda_{2,2}
+ut\lambda_{3,2} - z\lambda_{4,2} \\
= (z^2+&ut^2)(\alpha_{2,2}^2+\alpha_{1,2}^2u-2\alpha_{1,2}v+w) .
\end{aligned}
\end{equation}
Thus, if we localize at $z^2+ut^2\ne0$, we can add the following elements
to the ideal:
\begin{equation}
\begin{aligned}
\lambda_1 &= \alpha_{1,1}-\alpha_{2,2}\\
\lambda_2 &= 2v-\alpha_{1,2}u-\alpha_{2,1}\\
\lambda_3 &= \alpha_{2,2}^2+\alpha_{1,2}^2u-2\alpha_{1,2}v+w .
\end{aligned}
\end{equation}
Since the Grassmann blowup is irreducible, the elements of this extended
ideal must vanish on it, even when we do not impose $z^2+ut^2\ne0$.

In fact, the extended ideal is generated by 
$\lambda_{1,2}$, $\lambda_{2,2}$, $\lambda_1$, $\lambda_2$, and
$\lambda_3$, as we now verify: we have
\begin{equation}
\begin{aligned}
\lambda_{1,1}&=\lambda_{2,2}+z\lambda_1-t\lambda_2 \\
\lambda_{2,1}&=-u\lambda_{1,2}-ut\lambda_1-z\lambda_2 \\
\lambda_{3,1}&=(\alpha_{1,2}u-2v)\lambda_{1,2}+\alpha_{2,2}\lambda_{2,2}
+(x-vt)\lambda_1+y\lambda_2-z\lambda_3 \\
\lambda_{3,2}&=
-\alpha_{2,2}\lambda_{1,2}+\alpha_{1,2}\lambda_{2,2}+t\lambda_3\\
\lambda_{4,1}&=\alpha_{1,1}u\lambda_{1,2}+\alpha_{2,1}\lambda_{2,2}
+(-\alpha_{1,2}uz-\alpha_{2,2}ut+2vz)\lambda_1\\ &\quad+
(\alpha_{2,2}z-\alpha_{1,2}ut)\lambda_2-ut\lambda_3 
\\
\lambda_{4,2}&=\alpha_{1,2}u\lambda_{1,2}+\alpha_{2,2}\lambda_{2,2}
-z\lambda_3 .
\end{aligned}
\end{equation}
Moreover, since $\lambda_{1,2}=y+\alpha_{1,2}z+\alpha_{2,2}t$
and $\lambda_{2,2}=x+vt-\alpha_{1,2}ut+\alpha_{2,2}z$,
we can use the generators of the extended ideal to eliminate
$y$, $x$, $\alpha_{1,1}$, $\alpha_{2,1}$, and $w$, respectively, leaving 
$z$, $t$, $u$, $v$, $\alpha_{1,2}$ and $\alpha_{2,2}$ as coordinates
in this chart, with no further relations among them.

We can describe the map from the Grassmann blowup
back to the original hypersurface 
by setting each generator $\lambda_{1,2}$, $\lambda_{2,2}$, $\lambda_1$,
$\lambda_2$, $\lambda_3$ of the ideal to zero in turn, and solving the
resulting equation for an appropriate variable:
\begin{align}
\label{eq:y} y&=-\alpha_{1,2}z-\alpha_{2,2}t\\
\label{eq:x} x+vt&=\alpha_{1,2}ut-\alpha_{2,2}z\\
\label{eq:11} \alpha_{1,1} &= \alpha_{2,2}\\
\label{eq:21} \alpha_{2,1} &= 2v-\alpha_{1,2}u\\
\label{eq:w} w &= -\alpha_{2,2}^2-\alpha_{1,2}^2u+2\alpha_{1,2}v .
\end{align}

Equations \eqref{eq:y} and \eqref{eq:x} can be recast in matrix form:
\begin{equation}
\label{eq:matrixform}
\begin{bmatrix}
y\\x+vt
\end{bmatrix}
=
\begin{bmatrix}
-z & -t \\
ut & -z \\
\end{bmatrix}
\begin{bmatrix}
\alpha_{1,2}\\ \alpha_{2,2}\\
\end{bmatrix},
\end{equation}
and whenever $z^2+ut^2\ne0$, these
can be solved for $\alpha_{1,2}$ and
$\alpha_{2,2}$ in terms of $x$, $y$, $z$, $t$, $u$, $v$.
Equations \eqref{eq:11} and \eqref{eq:21} can then be used to solve
for $\alpha_{1,1}$ and $\alpha_{2,1}$, and when this has been
done, equation \eqref{eq:w} becomes equivalent to the original hypersurface
equation \eqref{eq:quad}.
  Thus, the fiber over a given
$(x,y,z,t,u,v,w)$ value is a single point unless $z^2+ut^2=0$.

When $z^2+ut^2=0$, if the coefficient
 matrix
\begin{equation} \label{alpha-matrix}
\begin{bmatrix}
-z & -t \\
ut & -z \\
\end{bmatrix}
\end{equation}
from equation \eqref{eq:matrixform}
has rank one, then the fiber over a given point takes the form
\begin{equation}
(\alpha_{1,2},\alpha_{2,2}) \mapsto (\alpha_{1,2}+ct, \alpha_{2,2}-cz).
\end{equation}
Substituting this into equation \eqref{eq:w} yields the additional relation
\begin{equation}
 \alpha_{2,2}z-\alpha_{1,2}ut+tv = 0
\end{equation}
which must hold if the fiber is nontrivial.  When this relation and
$z^2+ut^2=0$ both hold, the fiber is indeed one-dimensional.

On the other hand, if the matrix~\eqref{alpha-matrix} 
has rank $0$, then $z=t=0$
which implies $x=y=0$.  Any such point satifies 
equations \eqref{eq:y}--\eqref{eq:21}, leaving only equation \eqref{eq:w}, 
which describes---for
each fixed $(u,v,w)$---a
 conic embedded in the Grassmannian
mapping to $(0,0,0,0,u,v,w)$, as asserted at the
beginning of this section.

For later use, we summarize a portion of the structure of one of the
other coordinate charts, in less detail than above.
 For this second
chart,
we introduce four new variables
$\beta_{i,j}$ and eight new equations $\mu_{i,j}$ for the blowup by means of
\begin{equation}
\left[ \mu_{i,j} \right] = \Psi
\begin{bmatrix}
1 & 0 \\ \beta_{1,1} & \beta_{1,2} \\ 0 & 1 \\ \beta_{2,1} & \beta_{2,2} \\
\end{bmatrix} .
\end{equation}
Note that by multiplying by $\Phi$ from the left, we see that
the original equation is contained in the ideal generated by the
$\mu_{i,j}$.

We find some elements divisible by $y^2+wt^2=\partial W/\partial u$:
\begin{equation}
\begin{aligned}
y\mu_{1,1}+wt\mu_{1,2}+t\mu_{3,1}-y\mu_{3,2}
&=
(y^2+wt^2  )(\beta_{1,1}+\beta_{2,2}) ,
\\
-wt\mu_{1,1}+wy\mu_{1,2}+y\mu_{3,1}+wt\mu_{3,2}
&=
(y^2+wt^2  )(-\beta_{2,1}+\beta_{1,2}w) ,
\\
(-2vy+\beta_{1,2}wy+\beta_{2,2}wt)\mu_{1,2}
+(\beta_{1,2}&wt-\beta_{2,2}y)\mu_{3,2}
-wt\mu_{2,2}+y\mu_{4,2} \\
= (y^2+&wt^2)(\beta_{2,2}^2+u-2 \beta_{1,2}v+\beta_{1,2}^2w) .
\end{aligned}
\end{equation}
Localizing at $y^2+wt^2\ne0$, we can add the following elements to the
ideal:
\begin{equation}
\begin{aligned}
\mu_1 &= \beta_{1,1}+\beta_{2,2} \\
\mu_2 &= -\beta_{2,1}+\beta_{1,2}w \\
\mu_3 &= \beta_{2,2}^2+u-2 \beta_{1,2}v+\beta_{1,2}^2w .
\end{aligned}
\end{equation}
As before, the irreduciblity of the Grassmann blowup ensures that the
elements of the extended ideal vanish on it, even when we do not impose
$y^2+wt^2\ne0$.

In this chart, the extended ideal is generated by
$\mu_{1,2}$, $\mu_{3,2}$, $\mu_1$, $\mu_2$, and $\mu_3$, 
as we now demonstrate.
\begin{equation}
\begin{aligned}
\mu_{1,1}&=\mu_{3,2}+y\mu_1-t\mu_2\\
\mu_{2,1}&=(-2v+\beta_{1,2}w)\mu_{1,2}-\beta_{2,2}\mu_{3,2}+
(x+v t)\mu_1-z\mu_2-y\mu_3
\\
\mu_{2,2}&=\beta_{2,2}\mu_{1,2}+\beta_{1,2}\mu_{3,2}-t\mu_3\\
\mu_{3,1}&=-w\mu_{1,2}+wt\mu_1+y\mu_2\\
\mu_{4,1}&=\beta_{2,2}w\mu_{1,2}+\beta_{1,2}w\mu_{3,2}
-wz\mu_1-(x+v t)\mu_2-wt\mu_3\\
\mu_{4,2}&=(2v-\beta_{1,2}w)\mu_{1,2}+\beta_{2,2}\mu_{3,2}+y\mu_3 .
\end{aligned}
\end{equation}

Note that since
\begin{equation}
\begin{aligned}
 \mu_{1,2}&=\beta_{1,2}y+z+\beta_{2,2}t\\
 \mu_{3,2}&=\beta_{1,2}wt+x-vt-\beta_{2,2}y ,
\end{aligned}
\end{equation}
the Grassmann blowup
is also nonsingular in this coordinate chart: we can eliminate
$z$ with $\mu_{1,2}$, $x$ with $\mu_{3,2}$, $\beta_{1,1}$ with $\mu_1$,
$\beta_{2,1}$ with $\mu_2$, and $u$ with $\mu_3$.  This
leaves $y$, $t$, $v$, $w$, $\beta_{1,2}$, and $\beta_{2,2}$ as coordinates
in this chart.

\section{The $D_n$ case} \label{sec:6}

The main results of this paper concern
deformations of matrix factorizations for rational double
points of type $D_n$.  In this section, we 
will prove conjecture $3$ for the $D_n$ case,
from which conjectures $1$ and $2$ in the length $2$ case immediately
follow.  We will also prove the universality of the flop discussed
in section~\ref{sec:5} (i.e., Theorem~\ref{thm:universal}).

We follow notation for the versal deformation of $D_n$ and the invariant
theory for the corresponding Weyl group which was established in
\cite{gorenstein-weyl} (with one minor exception, noted below).
The matrix factorizations we use were first found in \cite{curto}.

The versal deformation of $D_n$ can be written
in the form 
\begin{equation}
X^2+Y^2Z-Z^{n-1} +2\gamma Y - \sum_{i=1}^{n-1}\delta_{2i}Z^{n-1-i}
= X^2+Y^2Z+2\gamma Y-F(Z),
\end{equation}
where $F(Z)$ is a general monic polynomial of degree $n-1$; the coefficients
of $F(Z)$, together with $\gamma$, give coordinates on the deformation
space $\operatorname{Def}$.

The invariant theory for $\mathfrak{W}_{D_n}$ can 
be described
in terms of the monic polynomial
\begin{equation}\label{originalpoly}
ZF(Z)+\gamma^2=
Z^n+\sum_{i=1}^{n-1}\delta_{2i}Z^{n-i} + \gamma^2=\prod_{j=1}^n(Z+t_j^2) ,
\end{equation}
which factors over the resolution space $\operatorname{Res}$ (the
coordinates on which are $t_1$, \dots, $t_n$).
The Weyl group $\mathfrak{W}_{D_n}$ is an extension of the symmetric group
$\mathfrak{S}_n$ (acting by permutations of the $t_j$'s) by a group 
$(\boldsymbol{\mu}_2)^{n-1}$.
The latter is the subgroup of $(\boldsymbol{\mu}_2)^n$ 
(acting on the $t_j$'s by 
coordinatewise multiplication) which preserves the product $t_1\dots t_n$.
The $\mathfrak{W}_{D_n}$-invariant functions are generated by
$\delta_{2i}=\sigma_i(t_1^2,\dots,t_n^2)$, the elementary symmetric functions
of $t_1^2$, \dots, $t_n^2$, together with\footnote{The sign in the definition
of $\gamma$ was not present in \cite{gorenstein-weyl}, but is convenient
here.  When comparing formulas in this paper to those in 
\cite{gorenstein-weyl}, one should replace $Y$ by $(-1)^nY$ to compensate
for this sign change.} 
$\gamma=(-1)^nt_1\cdots t_n$,
subject to
\begin{equation}
\delta_{2n}=\gamma^2 .
\end{equation}
This is the same set of functions given by
the coefficients of $F(Z)$ together with $\gamma$, verifying that
$\operatorname{Def}=\operatorname{Res}/\mathfrak{W}_{D_n}$ 
in this case.

Note that this setup continues to make
 sense for low values of $n$: $\mathfrak{W}_{D_3}$ is
an extension of $\mathfrak{S}_3$ by $(\boldsymbol{\mu}_2)^2$ which coincides
with $\mathfrak{W}_{A_3}=\mathfrak{S}_4$; $\mathfrak{W}_{D_2}$
is an extension of $\mathfrak{S}_2$ by $\boldsymbol{\mu}_2$ and coincides
with $\mathfrak{W}_{A_1\cup A_1}=\mathfrak{S}_2\times\mathfrak{S}_2$;
and $\mathfrak{W}_{D_1}$ is trivial.

The partial resolution corresponding to the $k^{\text{th}}$ vertex in
the Dynkin diagram has complementary graph $\Gamma_0$
of the form $\Gamma_{A_{k-1}}
\cup \Gamma_{D_{n-k}}$, for $1\le k\le n-2$, or $k=n$.  (The partial 
resolution for the $n{-}1^{\text{st}}$ vertex 
can be obtained from that for the $n^{\text{th}}$ by applying an 
automorphism to the
Dynkin diagram.)  The corresponding subgroup in the invariant theory is
\begin{equation}
\label{eq:Dsubgroup}
\mathfrak{S}_k \times \mathfrak{W}_{D_{n-k}}\subset\mathfrak{W}_{D_n},
\end{equation}
which can be considered for any $k$ between $1$ and $n$.  Note, however,
that this subgroup in the case of $k=n-1$ corresponds to the
complement of {\em two}\/ vertices in the Dynkin diagram---the two 
``short legs.''

The two polynomials whose coefficients (together with 
$\eta=(-1)^{n-k}t_{k+1}\cdots t_n$)
capture the invariant theory for the subgroup \eqref{eq:Dsubgroup}
are 
\begin{equation} \label{eq:twopolys}
 f(U)=\prod_{j=1}^k (U-t_j) \text{ and } 
Zh(Z)+\eta^2=\prod_{j=k+1}^n (Z+t_j^2) .
\end{equation}
In particular, the coefficients of $f$ and $h$, together with $\eta$,
give coordinates on the space 
$\operatorname{PRes} = \operatorname{Res}/\mathcal{W}(\Gamma_0)$.
To describe the map $\operatorname{PRes} \to \operatorname{Def}$,
we relate the coordinates on the two spaces, i.e., we relate the polynomials
in equation \eqref{eq:twopolys} to the original 
polynomial \eqref{originalpoly},
as follows.  First write
\begin{equation}
 f(U) = Q(-U^2) + U P(-U^2) ,
\end{equation}
encoding the coefficients of $f$ into the coefficients of $P$ and $Q$.
Note that if
$Z=-U^2$ then
\begin{equation}
\prod_{j=1}^k (Z+t_j^2) = f(U)f(-U) = Q(Z)^2 + Z P(Z)^2
\end{equation}
so that
\begin{equation}
\prod_{j=1}^n (Z+t_j^2) = (Zh(Z)+\eta^2)(Q(Z)^2 + Z P(Z)^2)
\end{equation}
and that $ \gamma = \eta Q(0)$, since $Q(0)=f(0)=(-1)^kt)1\cdots t_k$.
We also write $Q(Z)=ZS(Z)+Q(0)$ when needed: the coefficents of $S$, 
together with $Q(0)$, are equivalent to the coefficients of $Q$.

We can now put the polynomial 
\begin{equation}
 F(Z)=\frac1Z\left(-\gamma^2+\prod_{j=1}^n(Z+t_j^2)\right)
\end{equation}
which appears in
the versal deformation of $D_n$ into the form
\begin{equation} \label{eq:F}
 F(Z)=
h(Z)\big(Q(Z)^2+ZP(Z)^2\big) + \eta^2\big(
2Q(0)S(Z)+ZS(Z)^2 +P(Z)^2
\big) ;
\end{equation}
this, together with the formula $\gamma=\eta Q(0)$,
 provides the explicit map from $\operatorname{PRes}$ to
$\operatorname{Def}$.

For later use, we also
 observe that the relationship among $P(Z)$, $Q(Z)$, and $f(U)$
 is captured by the existence of
a polynomial in two variables $G(Z,U)$ satisfying
\begin{equation} \label{eq:Gdef}
UP(Z)+Q(Z)=(U^2+Z)G(Z,U)+f(U).
\end{equation}

This invariant theory analysis has provided precisely the functions we need in
order to
 use the universal length $2$ matrix factorization.
If we substitute
\begin{equation} \label{eq:mainsub}
\begin{gathered}
x=X,\ y=Y-\eta S(Z),\ z=Q(Z),\ t=P(Z),\\
 u=Z,\ v=\eta,\ w=-h(Z)
\end{gathered}
\end{equation}
in equation \eqref{eq:quad}, we find
\begin{equation}
\begin{gathered}
X^2{+}Z(Y{-}\eta S(Z))^2 {+} 2\eta(Y{-}\eta S(Z))Q(Z) {-} h(Z)Q(Z)^2
{-}(Zh(Z){+}\eta^2)P(Z)^2 \\
=X^2{+}Z(Y^2{-}\eta^2S(Z)^2) {+} 2\gamma(Y{-}\eta S(Z)) {-}h(Z)Q(Z)^2 
{-} (Zh(Z){+}\eta^2)P(Z)^2 \\
=X^2+Y^2Z+2\gamma Y-F(Z) ,
\end{gathered}
\end{equation}
using equation \eqref{eq:F}.
Thus, substituting \eqref{eq:mainsub} into \eqref{eq:mainfactor} gives
a matrix factorization of length $2$ defined over the space
$\operatorname{PRes}$.  This will turn out to be the matrix factorization
predicted by Conjecture~\ref{conj:3} for the $D_n$ case.

Note that over the origin of $\operatorname{PRes}$, our construction
specializes to a matrix 
factorization for the rational double point $D_n$, of the form
$\Phi = XI_4-\Xi$, $\Psi = XI_4+\Xi$, where
\begin{equation} \label{eq:matfacDneven}
\Xi = 
\begin{bmatrix}
0&Y&(-Z)^{\frac k2} &0\\
-YZ & 0 & 0 & (-Z)^{\frac k2} \\
(-1)^{n-1}(-Z)^{n-\frac{k+2}2} & 0 & 0 & -Y \\
0 & (-1)^{n-1}(-Z)^{n-\frac{k+2}2} & YZ & 0 \\
\end{bmatrix}
\end{equation}
if $k$ is even, and 
\begin{equation} \label{eq:matfacDnodd}
\Xi =
\begin{bmatrix}
0 & Y & 0 & (-Z)^{\frac{k-1}2} \\
-YZ & 0 & (-Z)^{\frac{k+1}2} & 0\\
0 & (-1)^{n-1}(-Z)^{n-\frac{k+3}2} & 0 & -Y \\
(-1)^{n-1}(-Z)^{n-\frac{k+1}2} & 0 & YZ & 0 \\
\end{bmatrix}
\end{equation}
if $k$ is odd.  These are the Gonzalez-Sprinberg--Verdier matrix
factorizations for $D_n$.

Note that the ``universal flop of length $2$'' from section~\ref{sec:5} is
a special case of our construction when $n=4$ and $k=2$.  In that case,
$f(U)=U^2+f_1U+f_0$ has degree $2$ and $h(Z)=Z+h(0)$ has degree $1$;
it follows that $Q(Z)=-Z+Q(0)$ and $P(Z)=P(0)$, where $P(0)=f_1$
and $Q(0)=f_0$.  We also have $S(Z)=-1$.  Thus,
\begin{equation}
\begin{gathered}
 x=X,\ y=Y+\eta,\ z=-Z+Q(0),\ t=P(0)\\
 u=Z,\ v=\eta,\ w=-Z-h(0) 
\end{gathered}
\end{equation}
or conversely,
\begin{equation}
\begin{gathered}
 X=x,\ Y=y-v,\ Z=u,\ \eta=v\\
 h(0)=-u-w,\ P(0)=t,\ Q(0)=z+u
\end{gathered}
\end{equation}
and so it is clear that $x$, $y$, $z$, $t$, $u$, $v$, $w$ is just another
set of coordinates for the space spanned by $X$, $Y$, $Z$, $\eta$, $h(0)$,
$P(0)$, $Q(0)$.  Thus, the universal flop of length $2$ and the $D_4$
deformation with $k=2$ coincide.

The classification theorem of \cite{gorenstein-weyl}, when specialized to
flops of length $2$, asserts that for every flop $Y\dashrightarrow Y^+$ 
of length $2$, after shrinking $Y$ there is a map $\pi: Y\to\Delta$ to the 
unit disk $\Delta$ and a map 
\begin{equation}
\rho: 
\Delta \to \operatorname{PRes}(\Gamma_{D_4} {-} \{\text{central vertex}\})
\end{equation}
such that the pullback via $\rho$ of the deformation of $D_4$ coincides
with $X$, and the pullback via $\rho$ of the universal partial resolution
coincides with $Y$.  Thus, since we have just shown that 
our universal flop of length $2$ coincides with
the deformation of $D_4$ over 
$\operatorname{PRes}(\Gamma_{D_4} {-} \{\text{central vertex}\})$,
Theorem~\ref{thm:universal} will follow once
we have verified that the Grassmann blowup in the case $k=2$, $n=4$
gives the simultaneous partial resolution, i.e., once we have verified
Conjecture~\ref{conj:3} in this case.

Returning to the case of general $n\ge4$, we next show that 
the maximal Cohen--Macaulay module associated to our matrix factorization
is reducible
for a few special
values of $k$ (for each fixed $n$). To see this, we introduce 
the following invertible change of basis matrices:
\begin{equation}
\begin{gathered}
B_0:=
\begin{bmatrix}
X-\eta&Y&Q(0)&1\\
-1&0&0&0\\
-Q(0)&1&0&0\\
Y&0&1&0\\
\end{bmatrix},
B_1:=
\begin{bmatrix}
1&0&0&0\\
-X-\eta&Y&Q(0)&1\\
-Q(0)&1&0&0\\
Y&0&1&0\\
\end{bmatrix},
\\
B_2 :=
\begin{bmatrix}
1&0&0&-\frac12\\
0&1&-\frac12Z&0\\
0&0&1&0\\
0&0&0&-1\\
\end{bmatrix},
B_3 :=
\begin{bmatrix}
1&0&1&0\\
0&-1&0&1\\
1&0&-1&0\\
0&-1&0&-1\\
\end{bmatrix}.
\end{gathered}
\end{equation}

When $k=1$, we have $P(Z)=1$, $S(Z)=0$ and $Q(Z)=Q(0)$, and we find
\begin{equation}
B_0 \Phi B_1^{-1} =
\begin{bmatrix}\varphi_1&0&0\\0&\varphi_2&0\\0&0&XI_2-\xi_1\end{bmatrix},
\text{ and }
B_1 \Psi B_0^{-1} = 
\begin{bmatrix}\psi_1&0&0\\0&\psi_2&0\\0&0&XI_2+\xi_1\end{bmatrix}
\end{equation}
where
$\varphi_2=\psi_1=1$,   
$\varphi_1=\psi_2=X^2 +Y^2Z +2\gamma Y-F(Z)$, and
\begin{equation}
\xi_1=
\begin{bmatrix}
\eta-Q(0)Y&-Z-Q(0)^2\\
Y^2-h(Z)&-\eta+Q(0)Y\\
\end{bmatrix}.
\end{equation}
Thus, the rank $2$
maximal Cohen--Macaulay module $(\Phi,\Psi)$ is a direct sum of
the rank $1$ module $(XI_2-\xi_1,XI_2+\xi_1)$ and
the trivial modules $(\varphi_1,\psi_1)$ and $(\varphi_2,\psi_2)$.
This is to be expected, since the corresponding vertex in
the Dynkin diagram has coefficient $1$ rather than $2$.

When $k=n$, we have $h(Z)=0$, $\eta=1$ and we find
\begin{equation}
B_2\Xi B_2^{-1} =
\begin{bmatrix}
\xi_2&0\\0&\xi_2\end{bmatrix}
\end{equation}
where
\begin{equation}
\xi_2=
\begin{bmatrix}
-P(Z)&Y-S(Z)\\
-YZ+ZS(Z)-2Q(Z)&P(Z) \\
\end{bmatrix}
\end{equation}
This time, the rank $2$ maximal Cohen--Macaulay module $(\Phi,\Psi)$ is
a direct sum of two copies of the rank $1$ module 
$(XI_2-\xi_2,XI_2+\xi_2)$.
Again, we expected a rank $1$ module due to the coefficient of this
vertex in the Dynkin diagram being $1$.

Finally, for $k=n-1$, we have $h(Z)=1$ and we find
\begin{equation}
B_3\Xi B_3^{-1} =
\begin{bmatrix}
\xi_3 & 0\\
0 & \xi_4\\
\end{bmatrix},
\end{equation}
where
\begin{equation}
\xi_3 =
\begin{bmatrix}
\eta P(Z)+Q(Z)& -Y+\eta S(Z)+P(Z)\\
Z(Y-\eta S(Z)+P(Z))+2\eta Q(Z) &\eta P(Z) -Q(Z)\\
\end{bmatrix}
\end{equation}
and
\begin{equation}
\xi_4 =
\begin{bmatrix}
-\eta P(Z)-Q(Z)& -Y+\eta S(Z)-P(Z)\\
Z(Y-\eta S(Z)-P(Z))+2\eta Q(Z) &\eta P(Z) +Q(Z)\\
\end{bmatrix}.
\end{equation}
This time there are two inequivalent summands of smaller rank, and indeed,
this subgroup of the Weyl group is {\em not}\/ obtained by deleting a single
vertex from the Dynkin diagram.  The two summands correspond to the two
short legs of the Dynkin diagram, and blowing up the given maximal
Cohen--Macaulay module will necessarily blow up both of the corresponding
curves.

We are now ready to carry out the Grassmann blowup, in the case of 
$k$ between $2$ and $n-2$, inclusive.
In the first chart, after using $\lambda_1$ and $\lambda_2$ to
eliminate $\alpha_{1,1}$ and $\alpha_{2,1}$, respectively, the 
Grassmann blowup is defined
by the ideal $(\lambda_{1,2},\lambda_{2,2},\lambda_3)$
where
\begin{equation}
\lambda_{1,2}=Y-\eta S(Z) + \alpha_{1,2}Q(Z)+\alpha_{2,2}P(Z)
\end{equation}
allows us to eliminate $Y$, and
\begin{equation}
\lambda_{2,2}=X + \eta P(Z)-\alpha_{1,2}ZP(Z) + \alpha_{2,2}Q(Z) 
\end{equation}
allows us to eliminate $X$.  The remaining generator is
\begin{equation}
\lambda_3=\alpha_{2,2}^2+\alpha_{1,2}^2Z-2\eta \alpha_{1,2}-h(Z) ,
\end{equation}
which defines the versal deformation of a $D_{n-k}$ 
singularity,\footnote{Note that we are including here the cases of
$D_3 = A_3$ and $D_2 = A_1 \cup A_1$, and that the deformation theory
works correctly in these ``degenerate'' cases.} 
since $h(Z)$ is a monic polynomial of degree $n-k-1$.

In the other chart, after using $\mu_1$ and $\mu_2$ to eliminate 
$\beta_{1,1}$ and $\beta_{2,1}$, respectively, the 
Grassmann blowup is defined
by the ideal $(\mu_{1,2},\mu_{3,2},\mu_3)$.
We have
\begin{equation}
\mu_{3,2}=X-\eta P(Z)- \beta_{1,2}h(Z)P(Z)-\beta_{2,2}(Y-\eta S(Z))
\end{equation}
which allows us to eliminate $X$ on this chart, but
\begin{equation}
\mu_{1,2}=Q(Z) + \beta_{1,2}(Y-\eta S(Z))+\beta_{2,2}P(Z)
\end{equation}
does not immediately allow elimination.  The third generator is
\begin{equation} \label{eq:mu3}
\mu_3=\beta_{2,2}^2+Z-2\eta \beta_{1,2}-\beta_{1,2}^2h(Z) .
\end{equation}

To understand the geometry of this chart, we follow an algebraic
trick introduced by Tyurina \cite{[T]}\footnote{This is a variant
of the original trick which Tyurina used to show the 
very existence of simultaneous 
resolutions for deformations of $D_n$ singularities.} and form the combination
\begin{equation}
\widetilde{\mu}_{1,2}=\mu_{1,2}-G(Z,\beta_{2,2})\mu_3 .
\end{equation}
where $G$ is the polynomial from equation \eqref{eq:Gdef}.
We then use the defining property for $G$ to compute:
\begin{equation}
\begin{aligned}
\widetilde{\mu}_{1,2}
&=(Y-\eta S(Z))\beta_{1,2} +f(\beta_{2,2}) 
+2\eta G(Z,\beta_{2,2})\beta_{1,2} +G(Z,\beta_{2,2})h(Z)\beta_{1,2}^2
\\
&
=(Y-\eta S(Z) + G(Z,\beta_{2,2})h(Z)\beta_{1,2} 
+ 2\eta G(Z,\beta_{2,2}))\beta_{1,2}
+f(\beta_{2,2}) .
\end{aligned}
\end{equation}
Thus, introducing the variable
\begin{equation}
\widetilde{Y}=Y-\eta S(Z) + G(Z,\beta_{2,2})h(Z)\beta_{1,2} 
+ 2\eta G(Z,\beta_{2,2}) ,
\end{equation}
we see that $\widetilde{\mu}_{1,2}$ defines
 a versal deformation of an $A_{k-1}$ singularity:
\begin{equation}
\widetilde{\mu}_{1,2}=\widetilde{Y}\beta_{1,2}+f(\beta_{2,2})
\end{equation}
since $f(U)$ is monic of degree $k$.  Note that $Z$ can be implicitly
eliminated using $\mu_3$ (i.e., equation \eqref{eq:mu3}).

Thus, the Grassmann blowup of $\Psi$ yields a space with a rational curve lying
over the origin, such that the central fiber has both an $A_{k-1}$ and
a $D_{n-k}$ singularity along that curve, and the deformation induces
versal deformations of these $A_{k-1}$ and $D_{n-k}$ singularities---in
fact, the same versal deformations encoded by our invariant theory
analysis.  A precisely analogous thing happens if we blow up
$\Phi$ instead of $\Psi$.  (This is clear, since $x\mapsto -x$ 
exchanges the two.)
This proves Conjecture~\ref{conj:3} in the $D_n$ case, and hence
also proves Conjectures~\ref{conj:1} and \ref{conj:2} for length $2$,
as well as Theorem~\ref{thm:universal}.

\section*{Appendix: Matrix Factorizations for $E_6$, $E_7$, and $E_8$}

\setlength{\unitlength}{1 true in}

In this appendix, we present the Gonzalez-Sprinberg--Verdier matrix
factorizations for $E_6$, $E_7$, and $E_8$.  Gonzalez-Sprinberg and Verdier
left some of the entries in the matrices undetermined, and we have
made choices for these.  We have also endeavored to make our matrices
agree with those found in Chapter 9 of \cite{yoshino}, after substituting
$Y=x$, $Z=y$, permuting the rows and columns, and subjecting the matrices
$\varphi_\ell^\bullet$ 
(defined below) to an overall sign change.

We write the equation for the rational double point in the form
\begin{equation}
f(X,Y,Z)=X^2+g(Y,Z),
\end{equation}
 and consider a maximal Cohen--Macaulay module of
length $\ell$.  Many such modules can be described in terms of a
matrix factorization of $\ell\times\ell$ matrices
over the ring $\mathbb{C}[Y,Z]$
\begin{equation} 
\varphi_\ell^\bullet \psi_\ell^\bullet = -g(Y,Z)I_\ell ,
\end{equation}
where 
$\bullet$
is a label (possibly empty) that is used to distinguish among different
modules of rank $\ell$ when needed.  Out of these matrices, we construct
\begin{equation}
 \Xi_\ell^\bullet = 
\begin{bmatrix} 
0 & \varphi_\ell^\bullet \\ -\psi_\ell^\bullet & 0 
\end{bmatrix} ,
\end{equation}
which determines the matrix factorization
\begin{equation}
(XI_{2\ell}-\Xi_\ell^\bullet)(XI_{2\ell}+\Xi_\ell^\bullet) 
= (X^2+g(Y,Z))I_{2\ell}
\end{equation}
for the rational double point.
In a few cases, the smaller matrices $\varphi_\ell^\bullet$ and 
$\psi_\ell^\bullet$
do not exist, and we give
 $\Xi_\ell^\bullet$ directly.

For each rational double point,
we have included a Dynkin diagram with vertices labeled by $\ell$ 
or $\ell^\bullet$ to make the explicit McKay
correspondence clear.

The problem of computing matrix factorizations 
$(\varphi_\ell^\bullet, \psi_\ell^\bullet)$ for $g(Y,Z)$ is known as
the computation of the {\em Auslander--Reiten quiver} \cite{AR}, and it
is in that context that the computations in Chapter 9 of \cite{yoshino} are
presented.

\subsection*{%
Case $E_6$: $g(Y,Z)=Y^3+Z^4$}
\quad

\medskip

\noindent
{\em Dynkin diagram:}
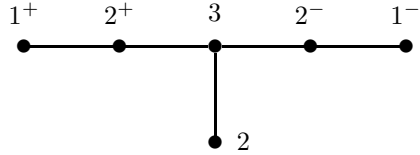
\begin{center}
\begin{picture}(2,1)(1.9,.5)
\thicklines
\put(1.9,1){\circle*{.075}}
\put(1.9,1){\line(1,0){.5}}
\put(2.4,1){\circle*{.075}}
\put(2.4,1){\line(1,0){.4625}}
\put(2.9,1){\circle*{.075}}
\put(2.9,.9625){\line(0,-1){.4625}}
\put(2.9,.5){\circle*{.075}}
\put(2.9375,1){\line(1,0){.4625}}
\put(3.4,1){\circle*{.075}}
\put(3.4,1){\line(1,0){.5}}
\put(3.9,1){\circle*{.075}}
\put(1.775,1.05){\makebox(.25,.25){\footnotesize $1^+$}}
\put(2.275,1.05){\makebox(.25,.25){\footnotesize $2^+$}}
\put(2.775,1.05){\makebox(.25,.25){\footnotesize $3$}}
\put(2.925,.375){\makebox(.25,.25){\footnotesize $2$}}
\put(3.275,1.05){\makebox(.25,.25){\footnotesize $2^-$}}
\put(3.775,1.05){\makebox(.25,.25){\footnotesize $1^-$}}
\end{picture}
\end{center}

\noindent
{\em Matrices:}
\[
\Xi_1^\pm = \begin{bmatrix}
\pm iZ^2 & -Y^2 \\
Y & \mp iZ^2 \\ 
\end{bmatrix}
\]
\[
\Xi_2^\pm = \begin{bmatrix}
\pm iZ^2 & 0 & -Y^2 & 0 \\
0 & \pm iZ^2 & YZ & -Y^2 \\
Y & 0 & \mp iZ^2 & 0 \\
Z & Y & 0 & \mp iZ^2 \\
\end{bmatrix}
\]
\[
\varphi_3 = \begin{bmatrix}
-Y^2 & -Z^3 & -YZ^2 \\ 
YZ & -Y^2 & Z^3 \\
Z^2 & -YZ & -Y^2 \\ 
\end{bmatrix}
\qquad
\psi_3 = \begin{bmatrix}
Y & 0 & -Z^2 \\ 
Z & Y & 0 \\ 
0 & -Z & Y \\
\end{bmatrix}
\]
\[
\varphi_2 = \begin{bmatrix} 
Y^2 & -Z^3\\
-Z & -Y\\
\end{bmatrix}
\qquad
\psi_2 = \begin{bmatrix} 
-Y & Z^3\\
Z & Y^2\\
\end{bmatrix}
\]

\subsection*{%
Case $E_7$: $g(Y,Z)=Y^3+YZ^3$}
\quad

\medskip

\noindent
{\em Dynkin diagram:}
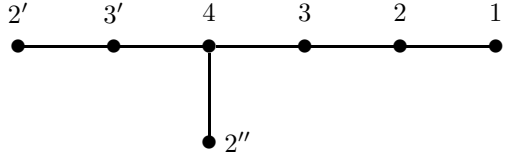
\begin{center}
\begin{picture}(3,1)(1.65,.5)
\thicklines
\put(1.9,1){\circle*{.075}}
\put(1.9,1){\line(1,0){.5}}
\put(2.4,1){\circle*{.075}}
\put(2.4,1){\line(1,0){.4625}}
\put(2.9,1){\circle*{.075}}
\put(2.9,.9625){\line(0,-1){.4625}}
\put(2.9,.5){\circle*{.075}}
\put(2.9375,1){\line(1,0){.4625}}
\put(3.4,1){\circle*{.075}}
\put(3.4,1){\line(1,0){.5}}
\put(3.9,1){\circle*{.075}}
\put(3.9,1){\line(1,0){.5}}
\put(4.4,1){\circle*{.075}}
\put(1.775,1.05){\makebox(.25,.25){\footnotesize $2'$}}
\put(2.275,1.05){\makebox(.25,.25){\footnotesize $3'$}}
\put(2.775,1.05){\makebox(.25,.25){\footnotesize $4$}}
\put(2.925,.375){\makebox(.25,.25){\footnotesize $2''$}}
\put(3.275,1.05){\makebox(.25,.25){\footnotesize $3$}}
\put(3.775,1.05){\makebox(.25,.25){\footnotesize $2$}}
\put(4.275,1.05){\makebox(.25,.25){\footnotesize $1$}}
\end{picture}
\end{center}

\noindent
{\em Matrices:}
\[
\varphi'_{2} = \begin{bmatrix} 
Y^2 & -YZ^2 \\
-Z & -Y\\ 
\end{bmatrix}
\qquad
\psi'_{2} = \begin{bmatrix} 
-Y & YZ^2 \\
Z & Y^2 \\ 
\end{bmatrix}
\]
\[
\varphi'_{3} = \begin{bmatrix}
Y^2 & -YZ^2 & Y^2Z \\ 
-YZ & -Y^2 & -YZ^2 \\
-Z^2 & -YZ & Y^2 \\ 
\end{bmatrix}
\qquad
\psi'_{3} = \begin{bmatrix}
-Y & 0 & YZ \\ 
Z & Y & 0 \\ 
0 & Z & -Y \\
\end{bmatrix}
\]
\[
\varphi_4 = \begin{bmatrix} 
0 & 0 & Y^2 & -YZ^2 \\
0 & 0 & -YZ & -Y^2 \\
Y & -Z^2 & 0 & YZ \\
-Z & -Y & -Y & 0 \\
\end{bmatrix}
\qquad
\psi_4 = \begin{bmatrix} 
0 & -YZ & -Y^2 & YZ^2 \\
Y & 0 & YZ & Y^2 \\
-Y & Z^2 & 0 & 0 \\
Z & Y & 0 & 0 \\
\end{bmatrix}
\]
\[
\varphi_3 = \begin{bmatrix}
-YZ & -Y^2 & -YZ^2 \\
Z^2 & YZ & -Y^2 \\ 
-Y & Z^2 & -YZ \\ 
\end{bmatrix}
\qquad
\psi_3 = \begin{bmatrix}
0 & -YZ & Y^2 \\ 
Y & 0 & -YZ \\
Z & Y & 0 \\ 
\end{bmatrix}
\]
\[
\varphi_2 = \begin{bmatrix} 
-YZ & Y^2 \\
-Y & -Z^2 \\ 
\end{bmatrix}
\qquad
\psi_2 = \begin{bmatrix} 
Z^2 & Y^2 \\ 
-Y & YZ \\
\end{bmatrix}
\]
\[
\varphi_1=\begin{bmatrix} -Y^2-Z^3 \end{bmatrix}
\qquad \psi_1=\begin{bmatrix} Y \end{bmatrix}
\]
\[
\varphi''_{2} = \begin{bmatrix} 
Y^2 & -YZ^2 \\
-YZ & -Y^2 \\ 
\end{bmatrix}
\qquad
\psi''_{2} = \begin{bmatrix} 
-Y & Z^2 \\
Z & Y \\ 
\end{bmatrix}
\]

\subsection*{%
Case $E_8$: $g(Y,Z)=Y^3+Z^5$}
\quad

\medskip

\noindent
{\em Dynkin diagram:}
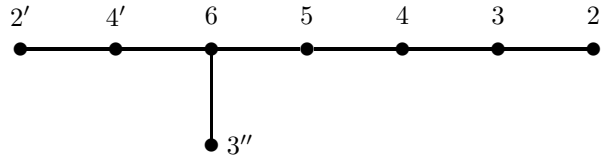
\begin{center}
\begin{picture}(3,1)(1.9,.5)
\thicklines
\put(1.9,1){\circle*{.075}}
\put(1.9,1){\line(1,0){.5}}
\put(2.4,1){\circle*{.075}}
\put(2.4,1){\line(1,0){.5}}
\put(2.9,1){\circle*{.075}}
\put(2.9,1){\line(0,-1){.5}}
\put(2.9,.5){\circle*{.075}}
\put(2.9,1){\line(1,0){.4625}}
\put(3.4,1){\circle*{.075}}
\put(3.4375,1){\line(1,0){.4625}}
\put(3.9,1){\circle*{.075}}
\put(3.9,1){\line(1,0){.5}}
\put(4.4,1){\circle*{.075}}
\put(4.4,1){\line(1,0){.5}}
\put(4.9,1){\circle*{.075}}
\put(1.775,1.05){\makebox(.25,.25){\footnotesize $2'$}}
\put(2.275,1.05){\makebox(.25,.25){\footnotesize $4'$}}
\put(2.775,1.05){\makebox(.25,.25){\footnotesize $6$}}
\put(2.925,.375){\makebox(.25,.25){\footnotesize $3''$}}
\put(3.275,1.05){\makebox(.25,.25){\footnotesize $5$}}
\put(3.775,1.05){\makebox(.25,.25){\footnotesize $4$}}
\put(4.275,1.05){\makebox(.25,.25){\footnotesize $3$}}
\put(4.775,1.05){\makebox(.25,.25){\footnotesize $2$}}
\end{picture}
\end{center}

\noindent
{\em Matrices:}
\[
\varphi'_{2} = \begin{bmatrix} 
-Z^3 & Y^2 \\
-Y & -Z^2 \\ 
\end{bmatrix}
\qquad
\psi'_{2} = \begin{bmatrix} 
Z^2 & Y^2 \\ 
-Y & Z^3 \\
\end{bmatrix}
\]
\[
\varphi'_{4} = \begin{bmatrix}
0 & -Z^3 & Y^2 & 0 \\
-Z^2 & 0 & -YZ & -Y^2 \\
Y & Z^2 & 0 & -Z^3 \\
0 & -Y & -Z^2 & 0 \\
\end{bmatrix}
\qquad
\psi'_{4} = \begin{bmatrix}
0 & Z^3 & -Y^2 & -YZ^2 \\
Z^2 & 0 & 0 & Y^2 \\
-Y & 0 & 0 & Z^3 \\
Z & Y & Z^2 & 0 \\
\end{bmatrix}
\]
\[
\varphi_6 = \begin{bmatrix}
0 & 0 & 0 & -Y^2 & -YZ^2 & -Z^4 \\
0 & 0 & 0 & -Z^3 & Y^2 & YZ^2 \\
0 & 0 & 0 & -YZ & -Z^3 & Y^2 \\
-Y & -Z^2 & 0 & 0 & 0 & -Z^3 \\
0 & Y & -Z^2 & Z^2 & 0 & 0 \\
-Z & 0 & Y & 0 & Z^2 & 0 \\
\end{bmatrix}
\]\[
\psi_6 = \begin{bmatrix}
0 & 0 & Z^3 & Y^2 & YZ^2 & Z^4 \\
-Z^2 & 0 & 0 & Z^3 & -Y^2 & -YZ^2 \\
0 & Z^2 & 0 & YZ & Z^3 & -Y^2 \\
Y & Z^2 & 0 & 0 & 0 & 0 \\
0 & -Y & Z^2 & 0 & 0 & 0 \\
Z & 0 & -Y & 0 & 0 & 0 \\
\end{bmatrix}
\]
\[
\varphi_5 = \begin{bmatrix}
Z^3 & Y^2 & 0 & 0 & 0 \\
0 & -Z^3 & -Y^2 & YZ^2 & Z^4 \\
0 & -YZ & Z^3 & Y^2 & YZ^2 \\
-Z^2 & 0 & -YZ & Z^3 & -Y^2 \\
Y & -Z^2 & 0 & 0 & 0 \\
\end{bmatrix}
\]\[
\psi_5 = \begin{bmatrix}
-Z^2 & 0 & 0 & 0 & -Y^2 \\
-Y & 0 & 0 & 0 & Z^3 \\
0 & Y & -Z^2 & 0 & 0 \\
-Z & 0 & -Y & -Z^2 & 0 \\
0 & -Z & 0 & Y & Z^2 \\
\end{bmatrix}
\]
\[
\varphi_{4} = \begin{bmatrix}
Z^3 & -Y^2 & 0 & 0 \\
0 & -YZ & Z^3 & Y^2 \\
Y & Z^2 & 0 & 0 \\
-Z & 0 & Y & -Z^2 \\
\end{bmatrix}
\qquad 
\psi_{4} = \begin{bmatrix}
-Z^2 & 0 & -Y^2 & 0 \\
Y & 0 & -Z^3 & 0 \\
0 & -Z^2 & -YZ & -Y^2 \\
Z & -Y & 0 & Z^3 \\
\end{bmatrix}
\]
\[
\varphi_3 = \begin{bmatrix}
-Y^2 & -Z^4 & -YZ^3 \\ 
-YZ & Y^2 & -Z^4 \\ 
-Z^2 & YZ & Y^2 \\
\end{bmatrix}
\qquad
\psi_3 = \begin{bmatrix}
Y & 0 & Z^3 \\
Z & -Y & 0 \\ 
0 & Z & -Y \\ 
\end{bmatrix}
\]
\[ 
\varphi_2 = \begin{bmatrix} 
Y^2 & -Z^4 \\
-Z & -Y \\ 
\end{bmatrix}
\qquad
\psi_2 = \begin{bmatrix} 
-Y & Z^4 \\
Z & Y^2 \\ 
\end{bmatrix}
\]
\[
\varphi''_{3} = \begin{bmatrix}
-Y^2 & -YZ^2 & -Z^4 \\ 
-Z^3 & Y^2 & YZ^2 \\
-YZ & -Z^3 & Y^2 \\ 
\end{bmatrix}
\qquad
\psi''_{3} = \begin{bmatrix}
Y & Z^2 & 0 \\
0 & -Y & Z^2 \\ 
Z & 0 & -Y \\ 
\end{bmatrix}
\]

\renewcommand{\MR}[1]{}
\bibliographystyle{amsplain-en}
\bibliography{flops}

\providecommand{\bysame}{\leavevmode\hbox to3em{\hrulefill}\thinspace}
\providecommand{\MR}{\relax\ifhmode\unskip\space\fi MR }
\providecommand{\MRhref}[2]{%
  \href{http://www.ams.org/mathscinet-getitem?mr=#1}{#2}
}
\providecommand{\href}[2]{#2}
\begin{thebibliography}{10}

\bibitem{Aspinwall}
P.~Aspinwall, \emph{A point's point of view of stringy geometry}, JHEP
  \textbf{01} (2003) 002, {\tt arXiv:hep-th/0203111}.

\bibitem{Aspinwall-Katz}
P.~Aspinwall and S.~Katz, \emph{Computation of superpotentials for {D}-branes},
  Comm. Math. Phys. \textbf{264} (2006) 227--253, {\tt arXiv:hep-th/0412209}.

\bibitem{Atiyah}
M.~F. Atiyah, \emph{On analytic surfaces with double points}, Proc. Roy. Soc.
  London. Ser. A \textbf{247} (1958) 237--244. \MR{MR0095974 (20 \#2472)}

\bibitem{AR}
M.~Auslander and I.~Reiten, \emph{Almost split sequences for
  {C}ohen--{M}acaulay-modules}, Math. Ann. \textbf{277} (1987) 345--349.
  \MR{MR886426 (88e:13002)}

\bibitem{bridgeland}
T.~Bridgeland, \emph{Flops and derived categories}, Invent. Math. \textbf{147}
  (2002) 613--632, {\tt arXiv:math.AG/0009053}. \MR{MR1893007 (2003h:14027)}

\bibitem{[Bri0]}
E.~Brieskorn, \emph{{\"U}ber die {A}ufl{\"o}sung gewisser {S}ingularit{\"a}ten
  von holomorphen {A}bbildungen}, Math. Ann. \textbf{166} (1966) 76--102.

\bibitem{CMII}
R.-O. Buchweitz, G.-M. Greuel, and F.-O. Schreyer, \emph{Cohen--{M}acaulay
  modules on hypersurface singularities. {II}}, Invent. Math. \textbf{88}
  (1987) 165--182. \MR{MR877011 (88d:14005)}

\bibitem{[CKM]}
H.~Clemens, J.~Koll{\'a}r, and S.~Mori, \emph{Higher dimensional complex
  geometry}, Ast\'erisque, vol. 166, Soci{\'e}t{\'e} Math{\'e}matique de
  France, Paris, 1988.

\bibitem{curto}
C.~Curto, \emph{Matrix model superpotentials and {C}alabi--{Y}au spaces: an
  {ADE} classification}, Ph.D. thesis, Duke University, 2005, {\tt
  arXiv:math.AG/0505111}.

\bibitem{Durfee}
A.~H. Durfee, \emph{Fifteen characterizations of rational double points and
  simple critical points}, Enseign. Math. (2) \textbf{25} (1979) 131--163.
  \MR{MR543555 (80m:14003)}

\bibitem{eisenbud}
D.~Eisenbud, \emph{Homological algebra on a complete intersection, with an
  application to group representations}, Trans. Amer. Math. Soc. \textbf{260}
  (1980) 35--64. \MR{MR570778 (82d:13013)}

\bibitem{GSV}
G.~Gonzalez-Sprinberg and J.~Verdier, \emph{Construction g\'eom\'etrique de la
  correspondence de {M}c{K}ay}, Ann. Sci. {\'E}cole Norm. Sup. (4) \textbf{16}
  (1983) 409--449.

\bibitem{Hori:2004zd}
K.~Hori and J.~Walcher, \emph{D-branes from matrix factorizations}, Comptes
  Rendus Physique \textbf{5} (2004) 1061--1070, {\tt arXiv:hep-th/0409204}.

\bibitem{Ito-Nakamura}
Y.~Ito and I.~Nakamura, \emph{Hilbert schemes and simple singularities},
  Euroconference (K.~Hulek et~al., eds.), Cambridge University Press, 1998,
  pp.~169--249.

\bibitem{Kapustin:2002bi}
A.~Kapustin and Y.~Li, \emph{D-branes in {L}andau--{G}inzburg models and
  algebraic geometry}, JHEP \textbf{12} (2003) 005, {\tt arXiv:hep-th/0210296}.

\bibitem{gorenstein-weyl}
S.~Katz and D.~R. Morrison, \emph{{G}orenstein threefold singularities with
  small resolutions via invariant theory for {W}eyl groups}, J. Algebraic Geom.
  \textbf{1} (1992) 449--530, {\tt arXiv:alg-geom/9202002}.

\bibitem{kawamata:crepant}
Y.~Kawamata, \emph{Crepant blowing-up of {$3$}-dimensional canonical
  singularities and its application to degenerations of surfaces}, Ann. of
  Math. (2) \textbf{127} (1988) 93--163. \MR{MR924674 (89d:14023)}

\bibitem{kawamata:flops}
Y.~Kawamata, \emph{General hyperplane sections of nonsingular flops in
  dimension $3$}, Math. Res. Lett. \textbf{1} (1994) 49--52, {\tt
  arXiv:alg-geom/9310002}.

\bibitem{CMI}
H.~Kn{\"o}rrer, \emph{Cohen--{M}acaulay modules on hypersurface singularities.
  {I}}, Invent. Math. \textbf{88} (1987) 153--164. \MR{MR877010 (88d:14004)}

\bibitem{[Kol]}
J.~Koll\'ar, \emph{Flops}, Nagoya Math. J. \textbf{113} (1989) 14--36.

\bibitem{kontsevich}
M.~Kontsevich, unpublished, as cited in \cite{Kapustin:2002bi}.

\bibitem{[L]}
H.~B. Laufer, \emph{On {$\mathbb{CP}^1$} as exceptional set}, Recent
  Developments in Several Complex Variables (J.~E. Fornaess, ed.), Ann. of
  Math. Stud., vol. 100, Princeton University Press, 1981, pp.~261--275.

\bibitem{McKay}
J.~McKay, \emph{Graphs, singularities, and finite groups}, {S}anta {C}ruz
  Conference on Finite Groups ({S}anta {C}ruz, 1979), Proc. Symp. Pure Math.,
  vol.~37, American Mathematical Society, 1980, pp.~183--186.

\bibitem{[P]}
H.~Pinkham, \emph{Factorization of birational maps in dimension 3},
  Singularities (P.~Orlik, ed.), Proc. Symp. Pure Math., vol. 40, part 2,
  American Mathematical Society, 1983, pp.~343--371.

\bibitem{pagoda}
M.~Reid, \emph{Minimal models of canonical 3-folds}, Algebraic Varieties and
  Analytic Varieties (S.~Iitaka, ed.), Adv. Stud. Pure Math., vol.~1,
  Kinokuniya, 1983, pp.~131--180.

\bibitem{[T]}
G.~N. Tyurina, \emph{Resolution of singularities of flat deformations of
  rational double points}, Functional Anal. Appl. \textbf{4} (1970) 68--73.

\bibitem{VdB}
M.~Van~den Bergh, \emph{Three-dimensional flops and noncommutative rings}, Duke
  Math. J. \textbf{122} (2004) 423--455, {\tt arXiv:math.AG/0207170}.
  \MR{MR2057015 (2005e:14023)}

\bibitem{yoshino}
Y.~Yoshino, \emph{Cohen--{M}acaulay modules over {C}ohen--{M}acaulay rings},
  London Mathematical Society Lecture Note Series, vol. 146, Cambridge
  University Press, Cambridge, 1990. \MR{MR1079937 (92b:13016)}

\end{thebibliography}

\end{document}